\documentclass[prl,aps,
twocolumn,
nofootinbib,
floatfix,
superscriptaddress]{revtex4-2}
\usepackage{graphicx}
\usepackage{epsfig}
\usepackage{rotating}
\usepackage{amssymb}
\usepackage{subfigure}
\usepackage{dsfont}
\usepackage{psfrag}
\usepackage{bbold}
\usepackage{xcolor}
\usepackage{multirow}
\usepackage{amsmath,euscript,array,mathrsfs}
\usepackage[utf8]{inputenc}
\usepackage[T1]{fontenc}
\usepackage{enumitem}
\setitemize{itemsep=4pt}

\usepackage[colorlinks=true,linktocpage=true,linkcolor=blue,citecolor=blue]{hyperref}

\usepackage{amsthm}
 
\newtheorem{theorem}{Theorem}

\newtheorem{corollary}{Corollary}[theorem]
\newtheorem{conjecture}[theorem]{Conjecture}

\newcommand{\be}{\begin{equation}}
\newcommand{\ee}{\end{equation}}
\newcommand{\bea}{\begin{eqnarray}}
\newcommand{\eea}{\end{eqnarray}}

\def\zmax{z_{\text{max}}}

\begin{document}

\title{Prime numbers and random walks in a square grid}

\author{Alberto Fraile}
\affiliation{Department of Control Engineering, Faculty of Electrical Engineering, Czech Technical University in Prague, Karlovo náměstí 13, 121 35, Czech Republic}
\email{albertofrailegarcia@gmail.com}
\author{Osame Kinouchi}
\affiliation{Universidade de São Paulo, Departamento de Física-FFCLRP, Ribeirão Preto, SP, Brazil}
\author{Prashant Dwivedi}
\affiliation{Department of Control Engineering, Faculty of Electrical Engineering, Czech Technical University in Prague, Karlovo náměstí 13, 121 35, Czech Republic}
\author{Roberto Martínez}
\affiliation{Euskal Herriko Unibertsitatea, Universidad del País Vasco, Barrio Sarriena, s/n 48940 Leioa, Spain}
\author{Theophanes E. Raptis}
\affiliation{Physical Chemistry Lab, Chemistry Dept., National Kapodistrian University of Athens, Greece}
\author{Daniel Fernández}
\affiliation{Science Institute, University of Iceland, Dunhaga 3, 107 Reykjavík, Iceland}
\email{licadnium@gmail.com}

%\date{\today}

\begin{abstract}
In recent years, computer simulations are playing a fundamental role in unveiling some of the most intriguing features of prime numbers. In this work, we define an algorithm for a deterministic walk through a two-dimensional grid that we refer to as Prime Walk (PW). The walk is constructed from a sequence of steps dictated by and dependent on the sequence of last digits of the primes. Despite the apparent randomness of this generating sequence, the resulting structure – both in 2d and 3d – created by the algorithm presents remarkable properties and regularities in its pattern that we proceed to analyze in detail.
\end{abstract}

\maketitle

%%%%%%%%%%%%%%%%%%%%%%%%%%%%%%%%%%%%%%%%
\section{Introduction}
\label{sec:Intro}

One can argue that prime numbers present perplexing features, in a hybrid of local unpredictability and global regular behavior. It is this interplay between randomness and regularity that motivated searches for both local and global patterns that could potentially become signatures for certain underlying fundamental mathematical properties. Patterns like the connections that are known to exist between the prime number sequence and the non-trivial zeros of the Riemann zeta function \cite{edwards}, are one of the most important open problems in mathematics \cite{luque}.

Since the formulation of the Riemann hypothesis, much has been done, yet much remains in the dark. It is often acknowledged that not a small number of mathematical discoveries have been accomplished after having assumed many conjectures or hypothesis to be valid a priori. For this reason, instead of attempting an analysis of the underlying behavior of the prime numbers, which has been an aspiration of mathematicians for centuries, we choose to perform ultra-large-scale computer calculations at a fundamental level.

In this work, presented as an unbiased computational experiment, we observe, rather than prescribe, exactly how the motion of a deterministic walk defined over the prime numbers sequence conspires to produce radically different results when compared to a simple random walk.

The algorithm defined below that creates our prime walk (PW) is simple, yet again, the PW itself appears to be complex and unpredictable. It can thus be placed as an example of emergent complexity. In fact, within the subfield of prime number theory, many examples \cite{stein, chmielewski, guariglia, cattani, vartziotis} can be found in which a simple algorithm defined over the prime number sequence give rise to complex structures spontaneously, a complexity within which regularities can be found, thus opening a new avenue for research into the distribution of prime numbers.

%%%%%%%%%%%%%%%%%%%%%%%%%%%%%%%%%%%%%%%%
\section{Methodology}
\label{sec:Methodology}

Here we propose what is to our knowledge an original way of number arrangement yielding an appealing visual structure in the form of a fractal plot. Inspired by Ulam’s spiral \cite{stein}, we assign positions to positive integers in a 2d plane following these rules:

\begin{itemize}[label=\textendash, leftmargin=0.3cm]
 \item The starting point is $(0, 0)$, assigned to $N=1$.
 \item Given the point $(x, y)$ assigned to number $N$, if $N+1$ is not a prime, the same point is assigned to it.
 \item If $N+1$ is a prime and its last digit is 1, we move up in the plane: $(x, y) \to (x, y+1)$.
 \item If $N+1$ is a prime and its last digit is 3, we move down in the plane: $(x, y) \to (x, y-1)$.
 \item If $N+1$ is a prime and its last digit is 7, we move to the left in the plane: $(x, y) \to (x-1, y)$.
 \item If $N+1$ is a prime and its last digit is 9, we move to the right in the plane: $(x, y) \to (x+1, y)$.
\end{itemize}

\noindent Note that the last digits of prime numbers are 1, 3, 7 and 9. The only exceptions are primes 2 and 5 at the very beginning of the algorithm. This can be easily implemented in a computer code. For details, see Appendix~\ref{appendix:a}. 

Of course, the choices above are arbitrary and could be modified, with a permutation of the moves for the different digits. However, it can be easily shown that any permutation necessarily leads to an equivalent result, the path described by the algorithm being a rotation or a mirror symmetry of the one resulting from the choice above. Equivalent algorithms can be ruled out.

%%%%%%%%%%%%%%%%%%%%%%%%%%%%%%%%%%%%%%%%
\section{Results}
\label{sec:Results}

\begin{figure}
	\centering
	\includegraphics[width=.45\textwidth]{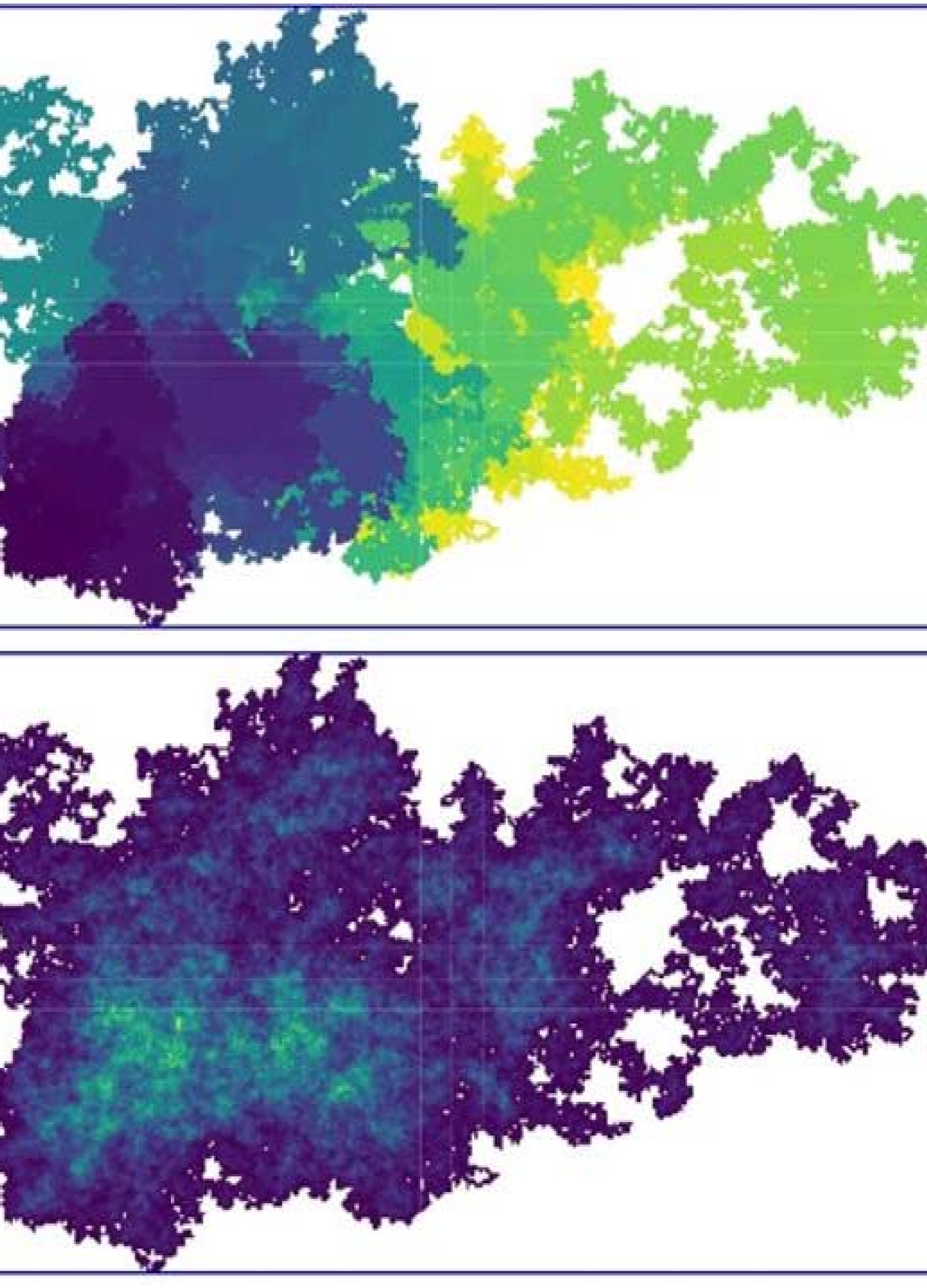}
	\vspace{0mm}
	\caption{(color online) Top: Example of the PW plot up to $N = 2 \cdot 10^7$. Color code represents the step. Bottom: As above, but now the color code represents the $\zmax$ value.}
	\label{Fig1}
\end{figure}

Following the algorithm, the walker will move through the grid in an erratic way, impossible to predict a priori. Fig.~\ref{Fig1} shows the PW created by the path up to $2 \cdot 10^7$. The color code is interpolated from the first steps (in dark blue) to the final steps (in yellow). Supplemental material to the published version of this article includes an animation showing how this area of the path grows with the increasing number of steps.

The PW is allowed to pass more than once through the same point in the grid. We can therefore keep track of the number of times that a certain point $(x, y)$ has been visited, and use this value as a third coordinate $z$ in order to visualize a structure in 3d. In supplemental material to the published version of this article, another animation is presented showing how the area covered by the path grows with the number of steps in terms of the maximum value of $z$ up to $N = 5 \cdot 10^{10}$. 

Furthermore, in order to help pinpoint patterns in our results, it can be interesting to compare them to those obtained from a random algorithm, in which at every prime $N$ the walker may move up, down, left or right in a random way (with equal probability). This produces a pseudo-random walk (pRW). Results from this alternative algorithm are presented in following figures along with the results obtained from the main prime walk.

%\vspace{0.5cm}

\begin{figure}
	\centering
	{\includegraphics[width=.47\textwidth ]{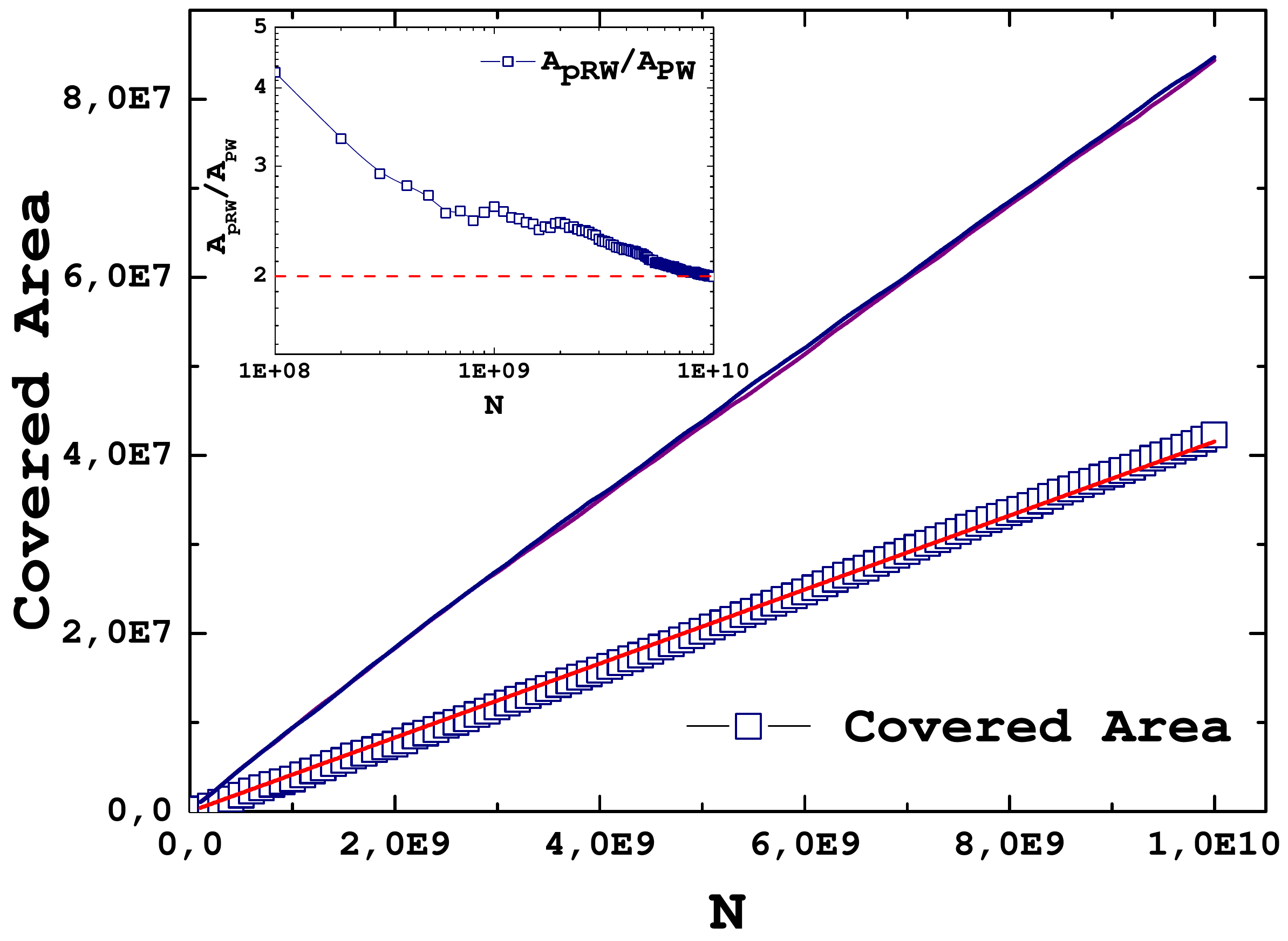}}
	\vspace{-4mm}
	\caption{Total area covered by the PW vs number of steps, $N$. We can see a linear trend (red line) is followed up to $10^{11}$ (Linear fit: $y = b x$). Purple and dark blue lines represent the result of two pseudo-random walks (pRW).}
	\label{Fig2}
\end{figure}

Finally, we also calculate the area covered by the path up to step $N$. Fig.~\ref{Fig2} plots this area vs the number of steps.
A linear scaling appears, with slope $b = 0.00414 \pm 9 \cdot 10^{-6}$, up to $N = 10^{11}$ steps. Whether this linear trend will hold for longer or start to saturate when a larger interval is explored is unclear.

\begin{figure}
	\centering
	{\includegraphics[width=.47\textwidth ]{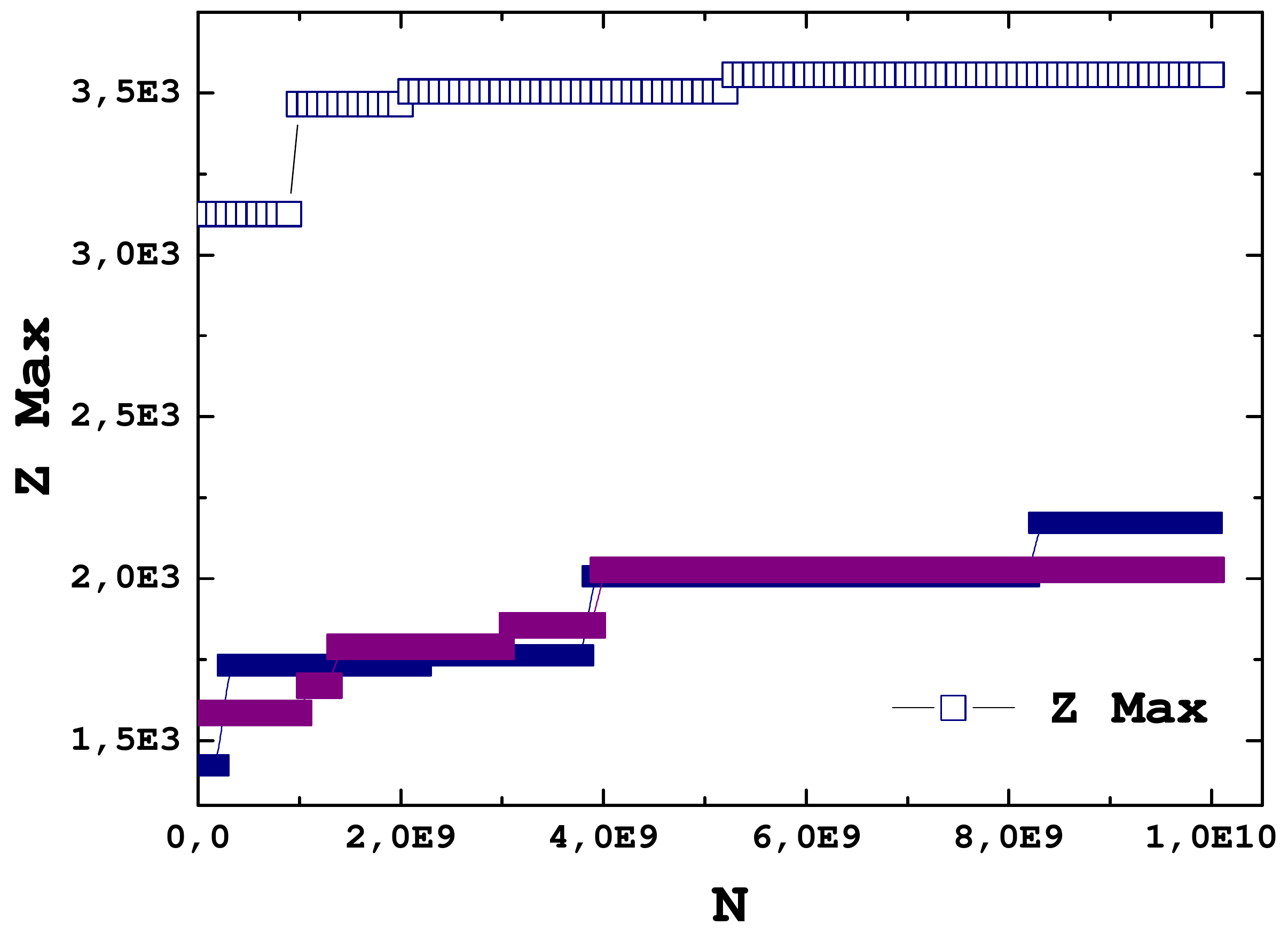}}
	\vspace{-3mm}
	\caption{Maximum $z$ value in the structure created by the PW vs number of steps, $N$. Full symbols are values corresponding to two pseudo-random walks.}
	\label{Fig3}
\end{figure}

The first thing one notices in Fig.~\ref{Fig3} is that, clearly, the area covered by the PW is smaller than the one covered by the pRWs. The difference is a \emph{factor of} 2 when $N$ is large enough (See inset).

It seems clear that the randomness of the prime number sequence produces an exact half-spread-path when the PW is compared to the pRW. This more compact RW is in contrast with the recently defined concept of maximum entropy random walk (MERW). As opposed to generic random walks (GRW), which maximize entropy locally (neighbors are chosen with equal probabilities), the MERW does it globally (all paths of given length and endpoints are equally probable) \cite{burda}.

The fact that, in a geometrical sense, the PW is more convoluted, spreading at a slower pace than the pRW, is also reflected in the maximum value of the $z$ coordinate, $\zmax$. This value can be computed, both in a cumulative way and within separate intervals. The differences between the values obtained from PW and pRWs are again clear in this case, being $\zmax$ higher for the PW by a factor that tends approximately to 1.6 (not shown).

Our computation is not large enough to produce a symmetric Gaussian distribution or Bell-shaped curve when we look at the structure constructed by the PW algorithm, however, it is natural to expect no preferences for any of the four quadrants of the plane.

%%%%%%%%%%%%%%%%%%%%%%%%%%%%%%%%%%%%%%%%
\section{Conjectures}
\label{sec:Conjectures}

Euclid's theorem tells us that there exists an infinite amount of primes, however, does this necessarily imply that the area covered by our PW path will turn out to be infinite?

\begin{figure}
	\centering
	{\includegraphics[width=.47\textwidth ]{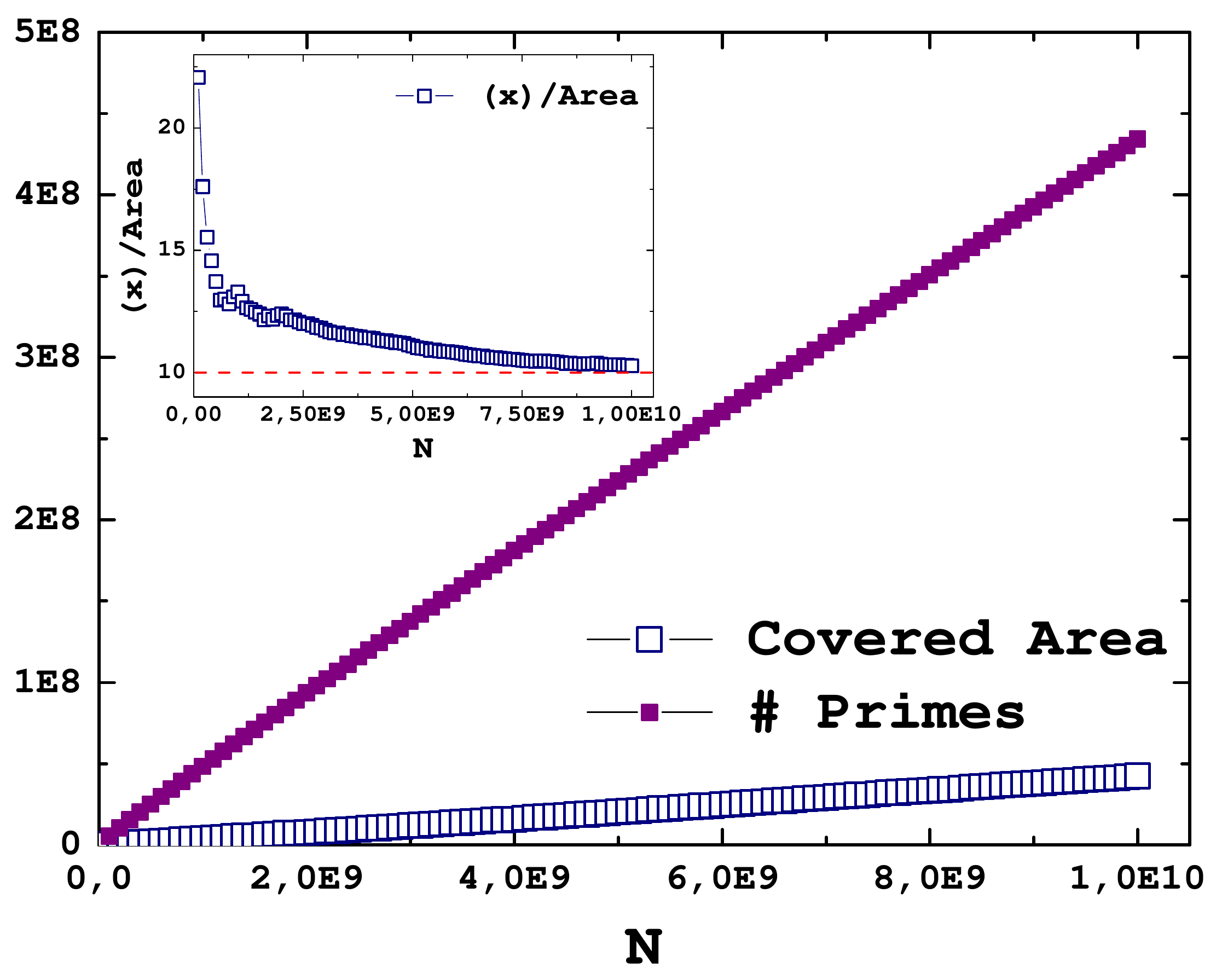}}
	\vspace{-3mm}
	\caption{(color online) Dark blue: area covered by the path up to $10^{10}$. Purple symbols: number of primes in the same interval, for comparison. Inset shows the ratio between both values vs $N$ - note the $\log \log$ scale.}
	\label{Fig4}
\end{figure}

It is natural to assume from the start that the area covered by the path after $N$ steps must be related in some way to the number of primes $\pi(N)$. In particular, we expect 
\be
\pi(N) \sim \frac{N}{\log N} \quad \text{when} \quad N \to \infty \,.
\ee
In our results, within the explored range, the covered area tends to a certain constant value $\psi$ times the number of primes, with this value being around $\psi = 1/10$ (see inset in Fig.~\ref{Fig4}). Apparently, the ratio
\be
\frac{\pi(N)}{A_{\text{pRW}}(N)} \to 10
\ee
when N is big enough. Will this be the case for even larger values of N?

We present here a plausible conjecture derived from our analysis, and two corollaries will follow straight forwardly. One of the main questions that our results invite us to ponder is whether there is an infinite number of points which are \emph{never} visited by the PW, or whether the whole 2d plane is visited at least once.

Consequently, we propose

\begin{conjecture}
\label{conj1}
The number of points within the area's perimeter which go unvisited by the PW increases homogeneously with $N$.
\end{conjecture}

\begin{corollary}
If conjecture \ref{conj1} is true, then there is an infinite number of points $(x, y)$ which are never visited by the PW. 
\end{corollary}
 
\begin{corollary}
If conjecture \ref{conj1} is true, then the area visited by the PW continues to grow indefinitely, ultimately becoming infinite in the limit $N \to \infty$.
\end{corollary}

In the explored range (up to $5 \cdot 10^{10}$) we observe an almost perfectly linear growth, but it seems clear that subsequently, with the primes becoming less frequent, the area will foreseeably continue growing at a lower rate. However, we conjecture that its growth will not stop (this is what Conjecture \ref{conj1} supports). The exact function describing this asymptotic growth is beyond the scope of this paper.

%%%%%%%%%%%%%%%%%%%%%%%%%%%%%%%%%%%%%%%%
\section{Discussion}
\label{sec:Discussion}

The covered area, as well as the value of $\zmax$, is governed partially (in a non-trivial way) by the gaps between primes. It is to be noted, though, that the "structure" created by the PW cannot be constrained in the $z$ direction by any upper limit, since the gaps between primes can be arbitrarily large. Nevertheless, this is just half of the problem, since for any given point $(x, y)$, it is impossible to know \emph{a priori} how many times it will have been revisited by the PW after a given number of steps $N$. Additionally, there is also the question of whether or not the PW is confined between some upper and lower values in $x$ and/or $y$.

Regarding the gaps between primes, it is known that gaps between consecutive prime numbers cluster on multiples of 6 \cite{wolf1, odlyzko}. Because of this, 6 is frequently called the ``jumping champion'', and it is conjectured that it holds this title all the way up to about 1035. Beyond 1035, and until 10425, the jumping champion then becomes 30 ($=2 \cdot 3 \cdot 5$), and beyond that, the most frequent gap is 210 ($=2 \cdot 3 \cdot 5 \cdot 7$) \cite{odlyzko}. Further important results on some statistical properties between gaps have been recently discovered \cite{szpiro1, szpiro2}. However, all of the aforementioned numerical observations, despite revealing intriguing properties about the prime sequence, cannot be easily applied to our problem to help figure out whether or not the PW will acquire definite boundaries.

On the other hand, according to Ares \emph{et al.} \cite{ares}, the apparent regularities previously observed \cite{wolf2, ball, kumar} reveal no structure in the sequence of primes; \emph{au contraire}, those regularities are precisely a consequence of its randomness. This is, however, a highly controversial topic. Recent computational work points out that ``after appropriate rescaling, the statistics of spacings between adjacent prime numbers follows the Poisson distribution''. See \cite{wolf3, garcia} and references therein for more on the statistics of the gaps between consecutive prime numbers.

\begin{figure}
	\centering
	{\includegraphics[width=.47\textwidth ]{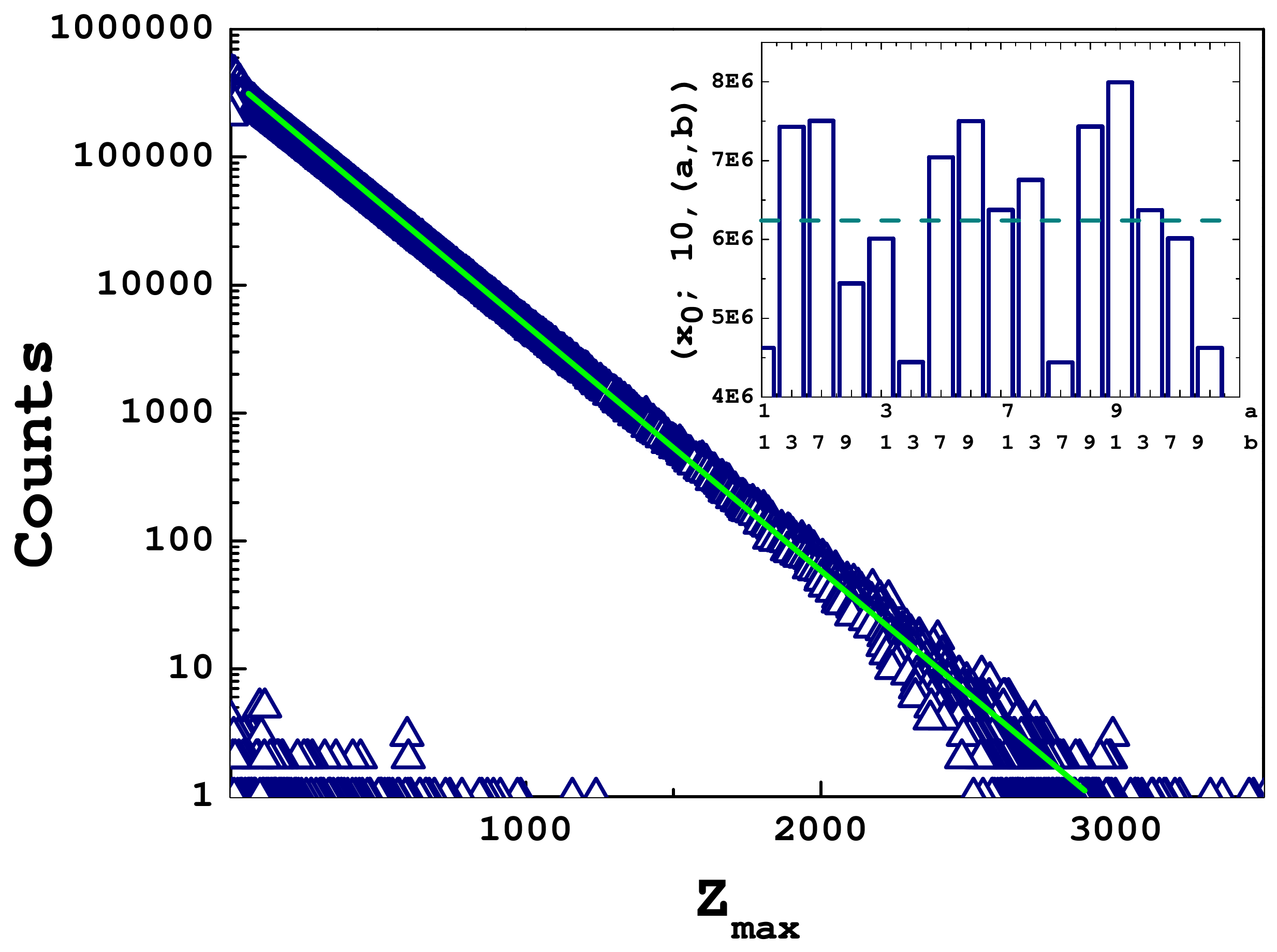}}
	\vspace{-2mm}
	\caption{(color online). Histogram of the values of $\zmax$ for $N = 10^{10}$. Note the $\log$ scale of the $x$ axis. Red line shows the linear fit of the data (see Inset). Among the first hundred million primes (modulo 10), there is substantial deviation from the prediction that each of the 16 pairs $(a, b)$ should have about $6.25 \cdot 10^6$ occurrences \cite{oliver}.}
	\label{Fig5}
\end{figure}

In Fig.~\ref{Fig5} we show the histogram $C(\zmax)$, that behaves almost as $\log C(\zmax) = b - a \zmax$ with  $a = - 0.0019$ and $b = 5.5$ for $N = 10^{10}$ (fit done removing low $x$, low $y$ points in the graph). This is similar to the histogram for differences between primes (prime gaps), see Fig. 1a of \cite{ares}. Note that $\zmax$ is a (non-random) sum of prime (but not necessarily consecutive) gaps. 

But we can be rather sure of a certain fact: these gaps, despite all of the complexity they present, lead to an absolutely clear Benford's law behavior for $\zmax$ (see Appendix~\ref{appendix:b} – Extra Analysis). This result cannot be coincidental. 

The second part of the problem concerns the last digit of the primes, a relatively unexplored topic. Many papers have been published about the first digits in the sequence of prime numbers \cite{luque, raimi, cohen}, but much less work has been devoted to studying the last digits. If we look at a pair of consecutive prime numbers $(a, b)$, assuming a purely random distribution one would expect it to be just as likely that these consecutive primes end in 1 and 1, or in 3 and 7, or in 3 and 9, and so on. However, intriguing irregularities (or biases) are actually observed in the distribution of consecutive primes. For example, it is a known fact that among the first hundred million primes (modulo 10), there is substantial deviation from the prediction that each of the 16 pairs $(a, b)$ should have about 6.25 million occurrences. Inset in Fig.~\ref{Fig5} shows the results according to Oliver and Soundararajan \cite{oliver}. Note that in our study we explored a little bit further since up to $10^{10}$ the number of primes is around $4.34 \cdot 10^8$.

This result must ultimately be the cause of the factor 2 observed in Fig.~\ref{Fig2}. The bias observed in \cite{oliver}, which essentially means that the prime numbers' last digits are usually not repeated when taken in consecutive pairs, contributes even more to the complexity of the PW when compared to the pRW. It actually does so in a very precise way, resulting in the covered area being exactly half the pRW value.

We know that the set of primes is algorithmically decidable since one can always find an effective primality test or a sieve to separate all primes. The problem of the same set having an underlying order or pattern is of a different nature. Specifically, such measures as Shannon's entropy may be insufficient in that they only count collective symbol occurrences. But we know that certain periodic as well as disordered patterns can give the same probabilities inside a whole multinomial set and hence, the exact same entropy value. For this reason a different measure of complexity as opposed to 'randomness' was proposed by Kolmogorov \cite{kolmogorov} in terms of the shortest formula or 'program' that can reproduce a given sequence. Although in principle an incalculable quantity, it can be approximated with data compression theory \cite{graham} which searches for redundancies this then being the equivalent of the shortest string able to reproduce the original. One can record in memory sufficiently large chunks of the last digit sequence and pass them through standard algorithms like Lempel-Ziv \cite{faloutsos} to find large compression ratios. Yet any such subset does not suffice for discovering a complete set of rules of fixed length for all primes.

%%%%%%%%%%%%%%%%%%%%%%%%%%%%%%%%%%%%%%%%
\section{Conclusions}
\label{sec:Conclusions}

In the present work we have intensively used a simple numerical representations of prime numbers in 2 and 3 dimensions to investigate the distribution of primes along the natural numbers.

Our mathematical experiment shows some important, unexpected and rather remarkable results. Within the explored range:

\begin{itemize}[leftmargin=0.4cm]
 \item The area covered by the PW is smaller than the one covered by the pRWs. Being the difference ``exactly'' a factor of 2 when $N$ is big enough.
 \item The number of primes up to $N$ is 10 times the area covered by the PW, $A(N)$. In other words, the covered area is 1/10 the number of primes $\pi(N)$.
 \item We showcase a remarkable match between the first digit count of the $\zmax$ values and Benford's law.
\end{itemize}

The results presented here highlight the important role of large-scale computer calculations as a way to discover \emph{possible} new properties of prime numbers. ``Possible'' needs to be stressed since we cannot prove that the results we observed will hold for any given larger range.

With the availability of increasing computational power, in a few years it will be possible to explore further orders of magnitude. This, however, will clearly never be enough. We need to turn our conjectures into demonstrated theorems.

It is interesting to note that the approach described both here and in our previous paper \cite{fraile} can be easily applied to any infinitely large sequence of numbers, such as for example the decimal digits of $\pi$, $e$, $\gamma$, or any other irrational number. Some interesting studies have been published about the randomness of $\pi$ \cite{marsaglia, aragon}, however, a plethora of questions remains open still. Could similar insights be extracted for them to some extent? We believe this to be an open question deserving of our attention.

\vspace{0.5cm}

%%%%%%%%%%%%%%%%%%%%%%%%%%%%%%%%%%%%%%%%
\begin{acknowledgments}
\noindent \emph{Acknowledgments}

\vspace{0.2cm}

This work was supported by The Ministry of Education, Youth and Sports from the Large Infrastructures for Research, Experimental Development and Innovations project ``IT4Innovations National Supercomputing Center – LM2015070''.

OK acknowledges CNAIPS-USP and CNPq. This work was produced as part of the activity of FAPESP Research, Innovation and Dissemination Center for Neuromathematics (grant \#2013/07699-0 S. Paulo Research Foundation).

\end{acknowledgments}

%%%%%%%%%%%%%%%%%%%%%%%%%%%%%%%%%%%%%%%%
%\begin{appendices}

\appendix
\section{Implementation}
\label{appendix:a}

\noindent \emph{PyPy implementation}

\vspace{0.2cm}

To satisfy the speed demand and reach the larger prime numbers, primality testing and prime numbers generating algorithms play a crucial role. The most crucial criteria in the analysis of the prime number generators are the number of probes, the number of generated primes, and the average time required in producing each prime. For our study we used the simple and efficient \emph{Sequences Containing Primes} algorithm, which employs function $m = 6k+1$ or $m = 6k-1$. The algorithm can be easily implemented in a code.

We used PyPy which is an implementation of the Python programming language written in RPython, a subset of the Python language, with its own interpreter \cite{bolz1, pypy, rpython}. It implements Python 2.7.10 and passes the Python test suite with some minor modifications \cite{python}. PyPy is intended to perform faster than CPython by employing a tracing Just-in Time (JIT) compiler. A JIT compiler, as the name suggests, compiles code during execution rather than before, as an Ahead-of-Time (AOT) compiler would do. The JIT used in PyPy is a meta-tracing JIT compiler. It does not encode any language semantics or profile in the execution of the program. Instead, it profiles the execution of the interpreter running the program. PyPy uses several optimizing techniques in their compiler; constant folding, common subexpression elimination, function inlining, and loop invariant code motion among others \cite{bolz2}. The trace also contains guards for each point in the recorded code that could branch off into another direction, for example in an if-statement. When the trace is compiled to machine code, each guard is compiled into a check that the execution is still correct. If it is not, the interpreter once again takes over execution. If a guard failure occurs more times than a certain limit, PyPy will attempt to compile the new execution branch as well.

\vspace{0.5cm}

\noindent \emph{Random implementation}

\vspace{0.2cm}

For the random implementation, we used probably the most widely known tool for generating random data in Python, its random module library. Python uses the Mersenne Twister PRNG algorithm \cite{matsumoto} as its core generator. Mersenne Twister (MT) was proposed for generating uniform pseudorandom numbers. For a particular choice of parameters, this algorithm provides a super astronomical period of $2^{19937} - 1$ and a 623-dimensional equidistribution up to 32-bit accuracy, while using a working area of only 624 words. The underlying implementation in C is both fast and threadsafe. MT is one of the most extensively tested random number generators in existence.

%%%%%%%%%%%%%%%%%%%%%%%%%%%%%%%%%%%%%%%%
\section{Extra Analysis}
\label{appendix:b}

\noindent \emph{Benford's Law}

\vspace{0.2cm}

\begin{figure}
	\centering
	{\includegraphics[width=.47\textwidth ]{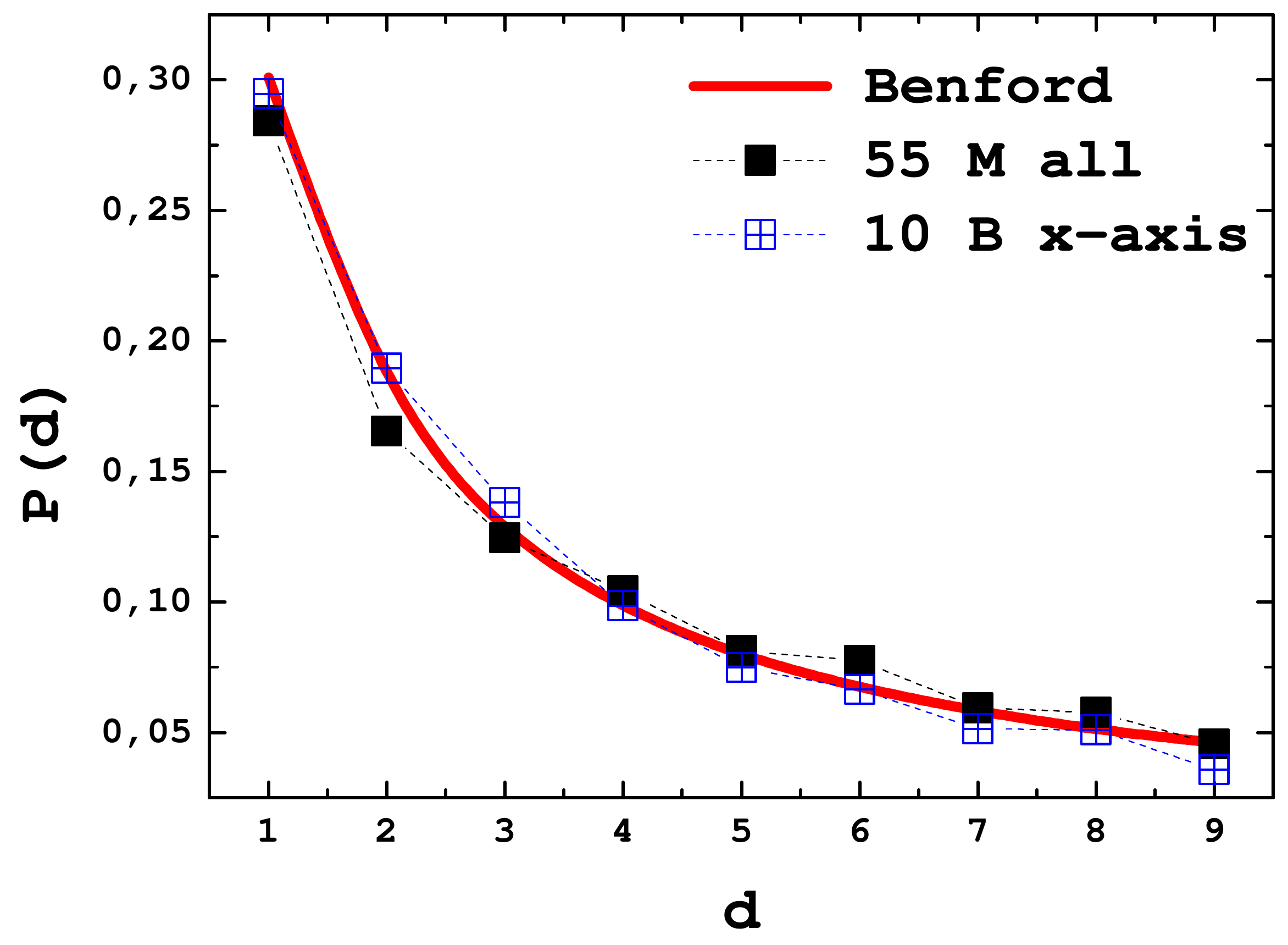}}
	\vspace{-3mm}
	\caption{(color online) Leading digit histogram of the $\zmax$ values (PW up to $5.5 \cdot 10^7$). The proportion of each of the $\zmax$ values is shown in black squares. Blue squares are another example, here up to $10^9$, but the $(x, y)$ points considered are only those in the $x$ axis. The expected values according to Benford's law are shown by the red curve.}
	\label{Figsm1}
\end{figure}

The first two figures in our manuscript invite us to ask: Are these z values randomly distributed, or do they possibly follow some kind of distribution? Could they actually follow Benford's law \cite{benfordslaw}? Fig.~\ref{Figsm1} below seems to indicate that to be the case (although a proof is beyond the scope of this paper) by plotting the $\zmax$ values even for small numbers from a statistical point of view.

For instance, when the PW reaches just $5.5 \cdot 10^7$, we have $\zmax$ = 155,802 points, and sorting these values according to the leading digit, comparing with Benford’s law the match is remarkable (Black squares in Fig.~\ref{Figsm1}). Even if we take just the $\zmax$ values along some given line (the $x$ axis, for example), the match is obvious (blue empty squares in Fig.~\ref{Figsm1}).

It is worth to note that the construct presented here resembles Jacob's Ladder \cite{fraile} (the description of the equivalence is beyond the scope of the present letter), so this result seems to point towards a behavior of the zeros in Jacob's Ladder following Benford's law as well when N is large enough.

\vspace{0.5cm}

\noindent \emph{Fractal dimension}

\vspace{0.2cm}

Fractal analysis is a contemporary method of applying non-traditional mathematics to describe patterns that defy understanding according to traditional Euclidean concepts. Recently, fractal analysis has been used to study a wide variety of complex patterns, like those of many types of biological cells \cite{kam}, tree and tumor growth \cite{cross}, gene expression \cite{aldrich}, forest fire progression \cite{turcotte}, economic trends, and cellular differentiation in space and time \cite{waliszewski}.

\begin{figure}
	\centering
	{\includegraphics[width=.47\textwidth ]{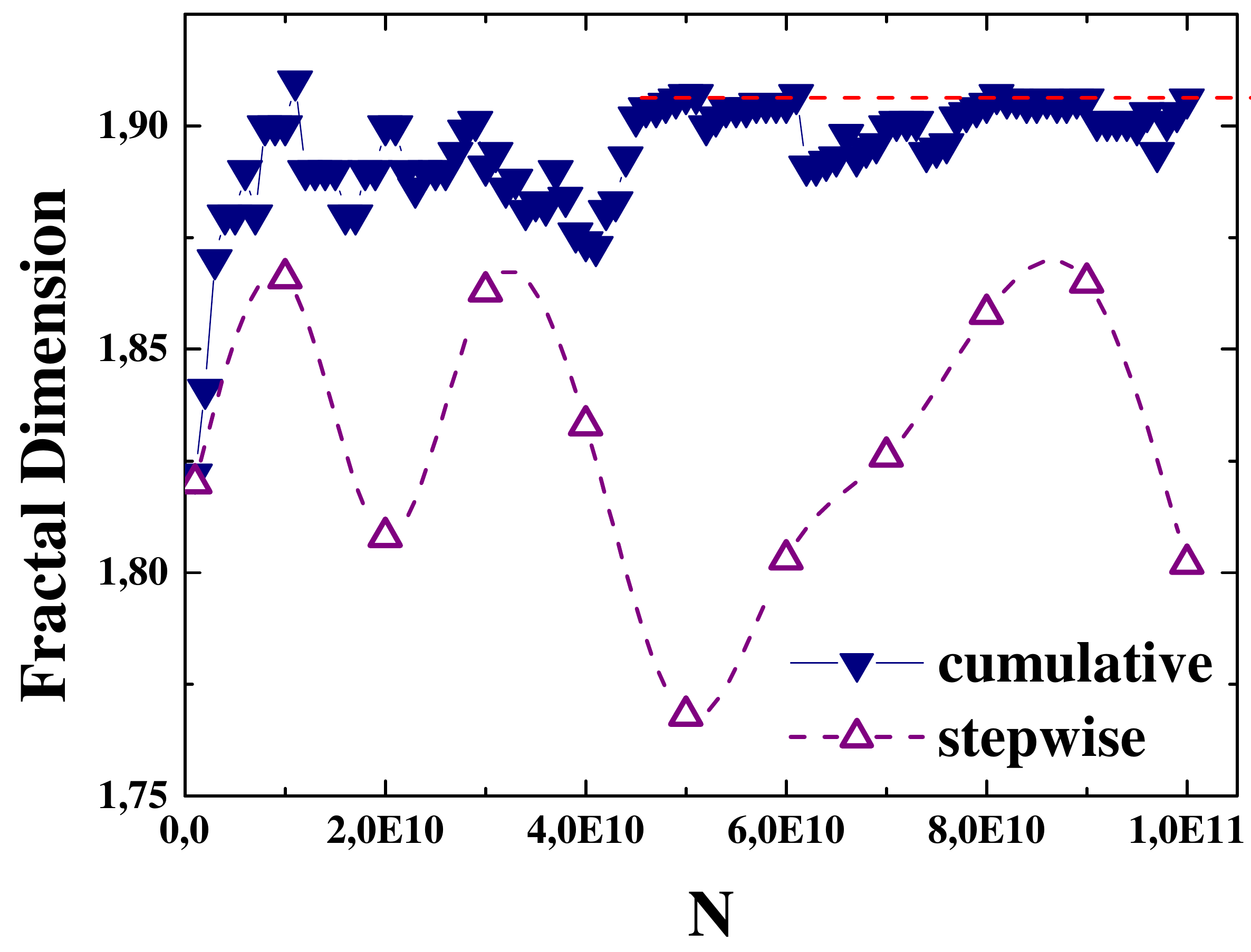}}
	\vspace{-3mm}
	\caption{Fractal dimension vs $N$ calculated in a cumulative way (blue inverted triangles) and in steps of $10^8$ (magenta empty triangles).}
	\label{Figsm2}
\end{figure}

In fractal analysis, complexity refers to the change in detail that comes with a change in scale. Many metrics of complexity can be defined, but the main parameter to capture them is the fractal dimension $D_F$ defined as a scaling rule comparing how detail in a pattern changes with the scale at which it is measured. Formally, each iteration driving the change in detail introduces new pieces into the fractal construct. The number of pieces $N$ at every step is related to the corresponding scale $\epsilon$ by
\be
N \propto \epsilon^{D_F}
\ee
Fig.~\ref{Figsm2} shows results for the fractal dimension $D_F$ as calculated with \emph{ImageJ}. Data points represent the result of the covered area vs $N$. We can safely conclude that $D_F$ tends to a value of 1.91 ($\pm 0.01$). On the other hand, we can see that if steps of $10^8$ are considered separately, larger oscillations are observed (as expected).

%\end{appendices}

%\newpage
%%%%%%%%%%%%%%%%%%%%%%%%%%%%%%%%%%%%%%%%

\end{document}